\documentclass[11pt]{amsart}

\usepackage{latexsym,amssymb,mathrsfs,amsmath,bbm,amsbsy,mathdots,color}
 \usepackage{hyperref}
\usepackage[english]{babel}
\usepackage[mathx]{mathabx}
\allowdisplaybreaks
\def\noteq{\operatorname{=\hspace{-3.2mm}/\ }}
 \newtheorem{theorem}{Theorem}[section]
 
 \newtheorem{lemma}[theorem]{Lemma}
 
 \theoremstyle{definition}
 
 \theoremstyle{remark}
 
 \numberwithin{equation}{section}

\begin{document}

\title[Paley-Wiener isomorphism]{Paley-Wiener isomorphism over\\ infinite-dimensional unitary groups}

\author[Lopushansky]{Oleh Lopushansky}

\address{%
1 Pigonia str.\\ 35-310 Rzesz\'ow\\
Poland}
\email{ovlopusz@ur.edu.pl; ovlopushansky@hotmail.com}

\thanks{Faculty of Mathematics and Natural Sciences, University of Rzesz\'ow}




\begin{abstract}
An analog of the Paley-Wiener isomorphism for the Hardy  space with an invariant
measure  over  infinite-dimensional unitary groups is described.
This  allows us to investigate  on such space the shift and multiplicative groups,
as well as, their generators and intertwining operators.
We show applications to the Gauss-Weierstrass semigroups and
to the Weyl-Schr{\"o}dinger irreducible representations
of complexified infinite-dimensional Heisenberg groups.
\end{abstract}

\maketitle
\tableofcontents

\section{Introduction}\label{Sec1}

The work deals with the Hardy space ${H}^2_\chi$ of  square-integrable complex-valued functions
with respect to a probability measure $\chi$ over the infinite-dimensional unitary group $U(\infty):=\bigcup\left\{U(m)\colon{m\in\mathbb{N}}\right\}$, extended by unit $\mathbbm{1}$,
which irreducibly acts on a separable complex Hilbert space ${E}$.
Here, $U(m)$ is the subgroup of unitary ${(m\times m)}$-matrices endowed with Haar's measure $\chi_m$.
In what follows, $U(\infty)$ is densely embedded  via a universal mapping $\pi$ into the space
of virtual unitary matrices $\mathfrak{U}=\varprojlim U(m)$ defined as the projective limit
under  Liv\v{s}ic's mappings $\pi_m^{m+1}\colon {U(m+1)\to U(m)}$.
The projective limit $\chi=\varprojlim\chi_m$,   such that each image-measure
$\pi_m^{m+1}(\chi_{m+1})$ is equal to $\chi_m$,  is concentrated on the range
$\pi(U(\infty))$ consisting of stabilized sequences
 (see  \cite[Neretin 2002]{zbMATH01820924}, \cite[Olshanski 2003]{zbMATH02037136}).
The measure  $\chi$ is invariant under right actions \cite[n.4]{zbMATH02037136}.
We refer to  \cite[Yamasaki 1974]{zbMATH03380964}, \cite[Borodin and Olshanski 2005]{zbMATH05004654}
for applications of $\chi$ to stochastic processes.
Needed properties of Hardy spaces ${H}^2_\chi$ can be found in \cite{lopushan2016}.  Various cases of Hardy spaces  in infinite-dimensional settings were considered in
\cite[Cole and Gamelin 1986]{ColeGamelin86}, \cite[{\O}rsted and Neeb 1998]{OrtedNeeb98}.

Now, we briefly describe results.
Using a  unitarily weighted symmetric Fock space $\varGamma_\mathsf{w}$,
defined by  ${E}$ and $\chi$, we find an orthogonal basis in  $H^2_\chi$ of Hilbert-Schmidt polynomials
such that   the conjugate-linear mapping
\[
\varPhi\colon\varGamma_\mathsf{w}\rightarrow{H}^2_\chi
\]
is a surjective isometry.  This allows us to establish in Theorem~\ref{laplace1}
an integral  formula for a Fock-symmetric $\mathcal{F}$-transform
\[
\mathcal{F}\colon{H}^2_\chi\ni f\mapsto\widehat{f}\in{H}^2_\mathsf{w}
\]
where the Hilbert space $H^2_\mathsf{w}$, uniquely determined by  $\varGamma_\mathsf{w}$,
consists of Hilbert-Schmidt analytic entire functions on $E$.
Thus, the  $\mathcal{F}$-transform acts as an analog of the Paley-Wiener isomorphism
over infinite-dimensional groups.

Furthermore, we investigate two different representations of the additive group from  $E$
over the Hardy space $H^2_\chi$  by shift and multiplicative groups.
Theorem~\ref{con} states that the $\mathcal{F}$-transform is an
intertwining operator between the multiplication group $M_\mathsf{a}^\dagger$ on ${H}^2_\chi$
and the shift group $T_\mathsf{a}$ on  ${H}^2_\mathsf{w}$.  On the other hand,
Theorem~\ref{conn} shows that $\mathcal{F}$ is the same between
the shift group $T_\mathsf{a}^\dagger$ on ${H}^2_\chi$ and
the multiplication group $M_{\mathsf{a}^*}$ on ${H}^2_\mathsf{w}$.
Integral formulas describing interrelations between their generators  are established.
In Theorem~\ref{comm} suitable commutation relations are stated.

Applications to the Gauss-Weierstrass-type semigroups on  ${H}^2_\chi$ are shown in Theorem~\ref{GW}.
An another application to linear representations of  complexified  infinite-dimensional Heisenberg groups on  ${H}^2_\chi$
in a Weyl-Schr\"{o}dinger form is given in Theorem~\ref{heis}.

Infinite-dimensional Heisenberg
groups was considered in \cite[Neeb 2000]{Neeb2000} by using reproducing kernel Hilbert spaces.
The Schr\"{o}dinger representation of infinite-dimensional Heisenberg groups
on $L^2_\gamma$ with respect to a Gaussian measure $\gamma$
over a real Hilbert space  is described in
\cite[I.Belti\c{t}\u{a}, D.Belti\c{t}\u{a} and  M.M\u{a}ntoiu 2016]{Beltita16}
(see also earlier publications  \cite{Beltita10,Beltita15}).

\section{Hilbert-Schmidt  analyticity}\label{Sec2}

Let $E$ stand for a separable complex Hilbert space with
scalar product ${\langle\cdot\mid\cdot\rangle}$ norm ${\|\cdot\|}$
and a fixed orthonormal basis $\left\{e_k\colon k\in\mathbb{N}\right\}$.
Denote by ${E}^{\otimes n}_{\sf alg}={E\otimes\stackrel{n \text{ times}}\ldots\otimes E}$ $(n\in\mathbb{N})$ its algebraic tensor power consisted of the linear span of elements $\psi_n={x_1\otimes\ldots\otimes x_n}$ with ${x_i\in{E}}$  ${(i=1,\ldots,n)}$. Set $x^{\otimes n}:={x\otimes\stackrel{n \text{ times}}\ldots\otimes x}$.
The  symmetric algebraic tensor power ${E}^{\odot n}_{\sf alg}={E\odot\ldots\odot E}$
is defined to be the range of the  projector
$\mathfrak{s}_n\colon{E}^{\otimes n}_{\sf alg}\ni\psi_n\mapsto
{x_1\odot\ldots\odot x_n}$ with ${x_1\odot\ldots\odot x_n} :=(n!)^{-1}
{\sum}_\sigma{x_{\sigma(1)}\otimes\ldots\otimes x_{\sigma(n)}}$ where
$\sigma\colon\{1,\ldots,n\}\mapsto\{\sigma(1),\ldots,\sigma(n)\}$
runs through all permutations. The symmetric algebraic Fock space is defined as the orthogonal sum
$\varGamma_{\sf alg}=
\bigoplus_{n\in\mathbb{Z}_+}{E}^{\odot n}_{\sf alg}$ with ${E}^{\odot 0}_{\sf alg}=\mathbb{C}$.

Let ${E}^{\otimes n}_\mathsf{h}:={E\otimes_\mathsf{h}\ldots\otimes_\mathsf{h} E}$
be the completion of ${E}^{\otimes n}_{\sf alg}$ by Hilbertian norm
${\left\|\psi_n\right\|_\mathsf{h}}={\langle\psi_n\mid\psi_n\rangle^{1/2}_\mathsf{h}}$
with  ${\left\langle\psi_n\mid\psi_n'\right\rangle_\mathsf{h}}=
{\langle x_1\mid x'_1\rangle\ldots\langle x_n\mid x'_n\rangle}$.
Denote by ${E}^{\odot n}_\mathsf{h}$
the range of continuous extension of $\mathfrak{s}_n$ on
${E}^{\otimes n}_\mathsf{h}$. As usual, the symmetric Fock space is defined to be
$\varGamma_\mathsf{h}=\bigoplus_{n\in\mathbb{Z}_+}{E}^{\odot n}_\mathsf{h}$.

Denote by $\lambda=(\lambda_1,\ldots,\lambda_m)\in\mathbb{N}^m$ with ${\lambda_1\ge\lambda_2\ge\ldots\ge\lambda_m}$
a partition of  ${n\in\mathbb{N}}$, that is, $n=|\lambda|$ where  $|\lambda|:={\lambda_1+\ldots+\lambda_m}$.
Any $\lambda$ may be identified with Young's diagram of length $l(\lambda)=m$.
Let $\mathbb{Y}$ denote all diagrams and $\mathbb{Y}_n=\left\{\lambda\in\mathbb{Y}\colon|\lambda|=n\right\}$.
Assume that  $\mathbb{Y}_0=\left\{\emptyset\in\mathbb{Y}\colon|\emptyset|=0\right\}$
and $l(\emptyset)=1$.
Let $\mathbb{N}^m_*:=
\left\{\imath=\left({\imath_1},\ldots,{\imath_m}\right)\in\mathbb{N}^m\colon\imath_l\noteq\imath_k,
\,\forall\,l\noteq k\right\}$. For each $\lambda\in\mathbb{Y}$ we assign the constant
\begin{equation}\label{constant}
C_{|\lambda|,l(\lambda)}:=\dfrac{(l(\lambda)-1)!|\lambda|!}{(l(\lambda)-1+|\lambda|)!}\le1.
\end{equation}

The spaces $E^{\odot n}_{\sf alg}$ and $\varGamma_{\sf alg}$ may be  generated by the basis  of symmetric tensors
\[\begin{split}
e^{\odot\mathbb{Y}_n}&=\bigcup\big\{e^{\odot\lambda}_\imath
:=e^{\otimes\lambda_1}_{\imath_1}\odot\ldots\odot{e}^{\otimes\lambda_{l(\lambda)}}_{\imath_{l(\lambda)}}\colon
(\lambda,\imath)\in\mathbb{Y}_n\times\mathbb{N}^{l(\lambda)}_\ast\big\}, \\
e^{\odot\mathbb{Y}}&=\bigcup\big\{e^{\odot\mathbb{Y}_n}\colon{n}\in\mathbb{Z}_+\big\}\quad\text{with}
\quad{e^{\odot\emptyset}_\imath=1},
  \end{split}\]
respectively. As is  known \cite[Sec. 2.2.2]{BerezanskiKondratiev95},
norm of basis element in  $\varGamma_\mathsf{h}$ is equal to
\begin{equation}\label{norm2fock}
\|{e}^{\odot\lambda}_\imath\|^2_\mathsf{h}=
\frac{\lambda!}{|\lambda|!},\qquad\lambda!:=\lambda_1!\cdot\ldots\cdot\lambda_m!.
\end{equation}

Let us define a new Hilbertian norm on $\varGamma_{\sf alg}$ by the equality $\|\cdot\|_\mathsf{w}=\langle\cdot\mid\cdot\rangle^{1/2}_\mathsf{w}$
where scalar product $\langle\cdot\mid\cdot\rangle_\mathsf{w}$ is determined via the
orthogonal relations
\begin{equation*}
\langle{e}^{\odot\lambda}_\imath\mid{e}^{\odot\lambda'}_{\imath'}\rangle_\mathsf{w}
=\left\{
\begin{array}{clcl}
 C_{|\lambda|,l(\lambda)}\|{e}^{\odot\lambda}_\imath\|^2_\mathsf{h} &: \lambda=\lambda&\text{ and }&\imath=\imath', \\
  0 &: \lambda\noteq\lambda'&\text{ or }&\imath\noteq\imath'.
\end{array}
\right.
\end{equation*}
Denote by ${E}^{\odot n}_\mathsf{w}$ and $\varGamma_\mathsf{w}$
the appropriate completions of $E^{\odot n}_{\sf alg}$  and $\varGamma_{\sf alg}$, respectively.
For any $\imath\in\mathbb{N}^{l(\lambda)}_\ast$  there corresponds
in ${E}^{\odot n}_\mathsf{w}$ the $d$-dimension  subspace  with $d=C_{|\lambda|,l(\lambda)}^{-1}$,
spanned by elements $\left\{e^{\odot\lambda}_\imath\colon\lambda\in\mathbb{Y}_n\right\}$.
The Hilbertian orthogonal  sum
\[
\varGamma_\mathsf{w}=\bigoplus_{n\in\mathbb{Z}_+}{E}^{\odot n}_\mathsf{w}
\]
endowed with $\langle\cdot\mid\cdot\rangle_\mathsf{w}$ we will call
{\em unitarily weighted} symmetric Fock space.

Let $x=\sum{e}_kx_k$ be the Fourier series of ${x\in E}$ with coefficients $x_k={\langle x\mid{e}_k\rangle}$.
We assign  to any $(\lambda,\imath)\in\mathbb{Y}_n\times\mathbb{N}^{l(\lambda)}_\ast$
the $n$-homogenous Hilbert-Schmidt polynomial defined via the Fourier coefficients
\[
x^\lambda_\imath:=\langle x^{\otimes n}\mid{e}^{\odot \lambda}_\imath\rangle_\mathsf{w}
=x^{\lambda_1}_{\imath_1}\ldots{x}^{\lambda_{l(\lambda)}}_{\imath_{l(\lambda)}},\quad
x\in{E}.
\]
Using the tensor multinomial theorem, we define  in $\varGamma_\mathsf{w}$ the Fourier decomposition
of exponential vectors (or coherent state vectors)
\begin{equation}\label{Tayl}
\begin{split}
\varepsilon(x)&:=\bigoplus_{n\in\mathbb{Z}_+}\frac{x^{\otimes n}}{n!}
=\bigoplus_{n\in\mathbb{Z}_+}\frac{1}{n!}\Big(\sum_{k\in\mathbb{N}}{e}_kx_k\Big)^{\otimes n}\\
&=\bigoplus_{n\in\mathbb{Z}_+}
\frac{1}{n!}\sum_{(\lambda,\imath)\in\mathbb{Y}_n\times\mathbb{N}^{l(\lambda)}_\ast}
\frac{n!}{\lambda !}{e}^{\odot\lambda}_\imath{x}^\lambda_\imath
\end{split}\end{equation}
with respect to  the basis ${e}^{\odot\mathbb{Y}}$. It is convergent in $\varGamma_\mathsf{w}$ in view of
 \eqref{constant} and
\begin{equation}\label{exp}\begin{split}
\|\varepsilon(x)\|^2_\mathsf{w}&=\sum_{n\in\mathbb{Z}_+}
\frac{1}{n!^2}\sum_{(\lambda,\imath)\in\mathbb{Y}_n\times\mathbb{N}^{l(\lambda)}_\ast}
\Big(\frac{n!}{\lambda !}\Big)^2\|{e}^{\odot\lambda}_\imath\|^2_\mathsf{w}|x^\lambda_\imath|^2\\
&=\sum_{n\in\mathbb{Z}_+}\frac{1}{n!^2}\sum_{(\lambda,\imath)}
\frac{n!}{\lambda !}C_{|\lambda|,l(\lambda)}|x^\lambda_\imath|^2
\le\sum_{n\in\mathbb{Z}_+}\frac{1}{n!}\sum_{(\lambda,\imath)}
\frac{n!}{\lambda !}|x^\lambda_\imath|^2\\
&=\sum_{n\in\mathbb{Z}_+}\frac{1}{n!}
\Big(\sum_{k\in\mathbb{N}}|x_k|^2\Big)^n=e^{\|x\|^2}.
\end{split}\end{equation}
Particulary, \eqref{exp} implies that the function ${E\ni x\mapsto\varepsilon(x)\in\varGamma_\mathsf{w}}$  is entire analytic.

Consider the space of complex-valued  functions in the variable ${x\in{E}}$
\[
H^2_\mathsf{w}:=\left\{\psi^*(x):=\left\langle\varepsilon(x)\mid\psi\right\rangle_\mathsf{w}\colon \psi\in\varGamma_\mathsf{w}\right\}
\quad\text{with norm}\quad\|\psi^*\|=\left\|\psi\right\|_\mathsf{w}.
\]
Every  function $\psi^*$ is entire analytic as the composition of $\varepsilon(\cdot)$ with
${\left\langle\cdot\mid\psi\right\rangle_\mathsf{w}}$.
The subspace in $H^2_\mathsf{w}$ of $n$-homogenous  Hilbert-Schmidt polynomials is defined to be
\[
H^{2,n}_\mathsf{w}=\big\{\psi_n^*(x)
=\langle x^{\otimes n}\mid\psi_n\rangle_\mathsf{w}\colon\psi_n\in{E}^{\odot n}_\mathsf{w}\big\}.
\]
Evidently, $H^2_\mathsf{w}=\mathbb{C}\oplus H^{2,1}_\mathsf{w}\oplus H^{2,2}_\mathsf{w}\oplus\ldots$.

It is important that $H^2_\mathsf{w}$
is uniquely determined by $\varGamma_\mathsf{w}$ since
$\left\{\varepsilon(x)\colon x\in{E}\right\}$  is total in $\varGamma_\mathsf{w}$.
Similarly, for the subspace $H^{2,n}_\mathsf{w}$ which is uniquely determined by ${E}^{\odot n}_\mathsf{w}$,
since $\left\{x^{\otimes n}\colon x\in{E}\right\}$  is total in ${E}^{\odot n}_\mathsf{w}$.
The last totality follows from the polarization formula  for symmetric tensor products
\begin{equation}\label{polar}
{e}^{\odot\lambda}_\imath=\frac{1}{2^nn!}
\sum_{\theta_ 1,\ldots,\theta_ n=\pm 1} \theta_1\dots \theta_n\,a^{\otimes n}\quad\text{with}\quad
a=\sum_{ i=1}^{l(\lambda)}\theta_ ie^{\otimes\lambda_i}_{\imath_i}
\end{equation}
which is valid for all ${{e}^{\odot\lambda}_\imath\in{e}^{\odot\mathbb{Y}_n}}$
(see e.g. \cite[Sec. 1.5]{zbMATH01445887})
Thus, the conjugate-linear isometries ${\psi\mapsto\psi^*}$ from $\varGamma_\mathsf{w}$ onto $H^2_\mathsf{w}$ and from ${E}^{\odot n}_\mathsf{w}$ onto $H^{2,n}_\mathsf{w}$ hold.

In conclusion, we can notice that every analytic  function ${\psi^*\in H^2_\mathsf{w}}$ determined by
$\psi={\bigoplus\psi_n\in\varGamma_\mathsf{w}}$,  $(\psi_n\in{E}^{\odot n}_\mathsf{w})$
has the Taylor expansion at zero
\[
\psi^*(x)=\sum_{n\in\mathbb{Z}_+}\frac{1}{n!}\sum_{(\lambda,\imath)\in\mathbb{Y}_n\times\mathbb{N}^{l(\lambda)}_\ast}
\frac{\langle{e}^{\odot\lambda}_\imath\mid\psi_n\rangle_\mathsf{w}}
{\|{e}^{\odot\lambda}_\imath\|^2_\mathsf{w}}
x^\lambda_\imath,\quad x\in E
\]
that follows from \eqref{Tayl}.
The function $\psi^*$ is entire  Hilbert-Schmidt analytic \cite[n.5]{lopushan2016}, \cite[n.2]{LopushanskyZagorodnyuk03}.

Note that analytic functions of Hilbert-Schmidt  types were considered in  \cite{Dwyer71}
 \cite{Petersson2001}.
More general classes of analytic functions associated with coherent sequences
of polynomial ideals were described in \cite{Carado09}.

\section{Hardy space over $U(\infty)$}\label{sec:3}

In what follows, we endow each  group $U(m)$ with the probability Haar measure $\chi_m$ and assume that
$U(m)$  is identified  with its range with respect to the embedding
$U(m)\ni u_m\mapsto  {\begin{bmatrix}
                      u_m & 0\\
                      0 &\mathbbm{1}\\
                    \end{bmatrix}\in U(\infty)}$.
The Liv\v{s}ic transform from $U(m+1)$ onto $U(m)$ is described in
\cite[Prop. 0.1]{zbMATH01820924} and \cite[Lem. 3.1]{zbMATH02037136} as the surjective Borel mapping
\begin{align*}\label{projective}
&\pi^{m+1}_m\colon{u_{m+1}}:=\begin{bmatrix}
                      z_m & a \\
                      b & t \\
                    \end{bmatrix}\longmapsto u_m:=\left\{\begin{array}{lc}
                      z_m-[a(1+t)^{-1}b]&: t\not\eq-1 \\
                      z_m&: t=-1. \\
                    \end{array}\right.
\end{align*}
The projective limit $\mathfrak{U}:=\varprojlim U(m)$
under $\pi^{m+1}_m$ has  surjective Borel projections $\pi_m\colon{\mathfrak{U}\ni u\mapsto u_m\in U(m)}$
such that $\pi_m={\pi^{m+1}_m\circ\pi_{m+1}}$.

Consider a universal dense embedding $\pi\colon U(\infty)\looparrowright\mathfrak{U}$
which to every ${u_m\in U(m)}$ assigns
the stabilized sequence $u=(u_k)$ such that  (see  \cite[n.4]{zbMATH02037136})
\begin{equation}\label{stab}
\pi\colon U(m)\ni u_m\mapsto(u_k)\in\mathfrak{U},\quad
u_k=\left\{\begin{array}{rl}
\pi^m_k(u_m)&:k<m\\
                      u_m &: k\ge m,
\end{array}\right.
\end{equation}
where $\pi^m_k:=\pi^{k+1}_k\circ\ldots\circ \pi^m_{m-1}$ for ${k<m}$
and $\pi^m_k$ is  identity mapping for ${k\ge m}$. On its range $\pi(U(\infty))$,
endowed with the Borel structure from $\mathfrak{U}$, we consider the  inverse mapping
 \[
 \pi^{-1}\colon\mathfrak{U}_\pi\rightarrow U(\infty)\quad
\text{where}\quad \mathfrak{U}_\pi:=\pi(U(\infty)).
\]
The right action $\mathfrak{U}_\pi\ni u\mapsto
u.g\in\mathfrak{U}_\pi$ with $g={(v,w)\in{U(\infty)\times U(\infty)}}$ is defined by
$\pi_m(u.g)=w^{-1}\pi_m(u)v$ where $m$ is so large that $g=(v,w)\in{U(m)\times U(m)}$.

Following \cite[n.3.1]{zbMATH01820924}, \cite[Lem. 4.8]{zbMATH02037136} via
 the Kolmogorov consistency theorem  (see e.g.  \cite[Thm 1]{Okada1977}, \cite[Cor. 4.2]{Tomas2006})
 we uniquely define on $\mathfrak{U}=\varprojlim U(m)$
 the probability measure $\chi:=\varprojlim\chi_m$   such that each image-measure
$\pi_m^{m+1}(\chi_{m+1})$ is equal to $\chi_m$.
For any Borel subset $A\subset\mathfrak{U}_\pi$ we have
$\pi_{m+1}(A)\subseteq(\pi^{m+1}_m)^{-1}\left[\pi_m(A)\right]$, because
$\pi_m={\pi^{m+1}_m\circ\pi_{m+1}}$. It follows that
${(\chi_m\circ\pi_m)(A)}=\pi^{m+1}_m(\chi_{m+1})[\pi_m(A)]=\chi_{m+1}[(\pi^{m+1}_m)^{-1}[\pi_m(A)]]
\ge{(\chi_{m+1}\circ\pi_{m+1})(A)}$. Hence,  $\chi$ satisfies the condition
\begin{equation}\label{occurs}
\chi(A)=\inf(\chi_m\circ\pi_m)(A)=\lim\chi_m(A)
\end{equation}
and therefore the projective limit
$\varprojlim\chi_m$ exists on $\mathfrak{U}_\pi$ via the well known Prohorov  theorem
\cite[Thm IX.52]{BourbakiINTII}.  Moreover, it is a Radon probability measure
concentrated on $\mathfrak{U}_\pi$  \cite[Thm 4.1]{Tomas2006}.
By the known portmanteau theorem \cite[Thm 13.16]{Klenke2008} and the Fubini theorem,
the invariance of Haar measures $\chi_m$ together with \eqref{occurs}
yield the  invariance properties under the right action,
\begin{align*}\label{inv0}
\int f(u.g)\,d\chi(u)&=\int f(u)\,d\chi(u),\quad g\in U(\infty)\times U(\infty),
\quad  f\in  L_\gamma^\infty,\\
\int f\,d\chi&=\int d\chi(\mathfrak{u})
\int_{U(m)\times U(m)}f(\mathfrak{u}.g)\,d(\chi_m\otimes\chi_m)(g),
\end{align*}
where $L_\chi^\infty$ stands for the space of all $\chi$-essentially bounded complex-valued functions
defined on $\mathfrak{U}_\pi$ and endowed with norm
$\|f\|_\infty=\mathop{\rm ess\,sup}_{u\in\mathfrak{U}_\pi}|f(u)|$.

Let $L^2_\chi$ be the space of square-integrable $\mathbb{C}$-valued functions $f$ on
$\mathfrak{U}_\pi$  with  norm
\[
\|f\|_\chi=\langle f\mid f\rangle_\chi^{1/2}\quad\text{where}\quad
\langle f\mid f\rangle_\chi:=\int f_1\bar f_2\,d\chi.
\]
The embedding $L^\infty_\chi\looparrowright L^2_\chi$  holds, moreover,
$\|f\|_\chi\le\|f\|_\infty$ for all $f\in{L}^\infty_\chi$.

To given the $E$-valued mapping $\mathfrak{U}_\pi\ni u\mapsto\pi^{-1}(u)e_1$,
we can well-define the Borel $\chi$-essentially bounded  functions in the variable ${u\in\mathfrak{U}_\pi}$,
\[
\phi_k:=\phi_{e_k},\quad
\phi_{e_k}(u)=\big\langle\pi^{-1}(u)e_1\mid e_k\big\rangle,
\quad k\in\mathbb{N},\]
which do not depend on the choice of $e_1$ in $\bigcup S(m)$
where $S(m)$ is the $m$-dimensional unit sphere in $E$ \cite[n.3]{lopushan2016}.
The uniqueness of $\phi_x(u)=\langle\pi^{-1}(u)e_1\mid{x}\rangle$ with $x\in E$ results from the total embedding $\pi\colon U(\infty)\looparrowright\mathfrak{U}$.
From \eqref{stab} it follows that $\pi^{-1}\circ\pi_m^{-1}$ coincides with the embedding
$U(m)\looparrowright U(\infty)$. Hence, by \eqref{occurs} and  the portmanteau theorem there exist the limit
\[
\int\phi_x\,d\chi=\lim_{m\to\infty}\int_{U(m)}\!\phi_x\,d(\chi_m\circ\pi_m)=
\lim_{m\to\infty}\int_{U(m)}\!(\phi_x\circ\pi_m^{-1})\,d\chi_m,
\]
i.e., $\phi_x\in L^\infty_\chi$ for any  $\phi_x(u)=\langle\pi^{-1}(u)e_1\mid{x}\rangle$ with
${x\in E}$.

By  formula  \eqref{polar}
to every ${e^{\odot\lambda}_\imath\in{e}^{\odot\mathbb{Y}_n}}$
there uniquely corresponds the Borel function from $L^\infty_\chi$
\begin{equation*}\label{base2}
\phi^\lambda_\imath(u):=\big\langle[\pi^{-1}(u)e_1]^{\otimes n}\mid
e^{\odot \lambda}_\imath\big\rangle_\mathsf{w}=
\phi^{\lambda_1}_{\imath_1}(u)\ldots\phi^{\lambda_{l(\lambda)}}_{\imath_{l(\lambda)}}(u)
\end{equation*}
in the variable ${u\in\mathfrak{U}_\pi}$.
It follows that the orthogonal basis $e^{\odot\mathbb{Y}}$ of elements ${e}^{\odot\lambda}_\imath=
{{e}^{\otimes\lambda_1}_{\imath_1}\odot\ldots\odot{e}^{\otimes\lambda_m}_{\imath_m}}$,
indexed by $\lambda={(\lambda_1,\ldots,\lambda_m)\in\mathbb{Y}}$ and
$\imath={\left({\imath_1},\ldots,{\imath_m}\right)\in\mathbb{N}^m_\ast}$ with $m=l(\lambda)$,
uniquely determines the systems
of Borel $\chi$-essentially bounded  functions in the variable  $u\in\mathfrak{U}_\pi$,
\[\begin{split}
\phi^{\mathbb{Y}_n}&=\bigcup\big\{\phi^\lambda_\imath:=
\phi^{\lambda_1}_{\imath_1}\cdot\ldots\cdot\phi^{\lambda_m}_{\imath_m}
\colon(\lambda,\imath)\in\mathbb{Y}_n\times\mathbb{N}^m_\ast,\ m=l(\lambda)\big\}, \\
\phi^\mathbb{Y}&=\bigcup\big\{\phi^{\mathbb{Y}_n}\colon{n}\in\mathbb{Z}_+\big\}\quad\text{with}
\quad{\phi^\emptyset_\imath\equiv1}.
  \end{split}\]

{\em The  Hardy space} $H_\chi^2$  is defined as the closed complex linear span of
$\phi^\mathbb{Y}$ endowed with $L^2_\chi$-norm. The following assertion is  proved in
 \cite[Thm 3.2]{lopushan2016}.

\begin{theorem}\label{irrep1}
The system of Borel functions $\phi^\mathbb{Y}$ forms an orthogonal basis in $H_\chi^2$  such that
\begin{equation*}\label{norm}
\|\phi^\lambda_\imath\|_\chi=C_{|\lambda|,l(\lambda)}^{1/2}
\|{e}^{\odot\lambda}_\imath\|_\mathsf{h},\quad
\lambda\in\mathbb{Y},\quad\imath\in\mathbb{N}^{l(\lambda)}_\ast.
\end{equation*}
\end{theorem}

Define the subspace $H_\chi^{2,n}\subset H_\chi^2$ for any ${n\in\mathbb{N}}$
to be the closed linear span of the subsystem $\phi^{\mathbb{Y}_n}$. Theorem~\ref{irrep1} implies  that
$H_\chi^{2,n}\perp H_\chi^{2,m}$  in $L^2_\chi$ for any  $n\noteq m$.  This provides the orthogonal decomposition
\[H_\chi^2=\mathbb{C}\oplus H_\chi^{2,1}\oplus H_\chi^{2,2}\oplus\ldots.\]

\section{Fock-symmetric $\mathcal{F}$-transform}\label{4}

The one-to-one correspondence $e^{\odot \lambda}_\imath\rightleftarrows\phi^\lambda_\imath$ allows us to define
via the change of orthonormal bases
\[\varPhi\colon\varGamma_\mathsf{w}\ni{e}^{\odot \lambda}_\imath\|e^{\odot\lambda}_\imath\|^{-1}_\mathsf{w}
\mapsto\phi^\lambda_\imath\|\phi^\lambda_\imath\|^{-1}_\chi\in{H}^2_\chi,
\qquad \lambda\in\mathbb{Y},\quad \imath\in\mathbb{N}^{l(\lambda)}_\ast
\]
the isometric conjugate-linear mapping $\varPhi\colon\varGamma_\mathsf{w}\rightarrow{H}^2_\chi$.
The adjoint  mapping ${\varPhi^*\colon{H}^2_\chi\rightarrow\varGamma_\mathsf{w}}$
is defined by  $\left\langle\varPhi e^{\odot\lambda}_\imath\mid f\right\rangle_\chi=
\left\langle{e}^{\odot\lambda}_\imath\mid\varPhi^*f\right\rangle_\mathsf{w}$ with ${f\in H^2_\chi}$.
The suitable Fourier decomposition has the form
\[
\varPhi\psi=\sum_{(\lambda,\imath)\in\mathbb{Y}\times\mathbb{N}^{l(\lambda)}_\ast}
\hat\psi_{(\lambda,\imath)}\phi^\lambda_\imath\|\phi^\lambda_\imath\|^{-1}_\chi,
\quad \hat\psi_{(\lambda,\imath)}:=\langle{e}^{\odot\lambda}_\imath\mid\psi\rangle_\mathsf{w}\,
\|{e}^{\odot\lambda}_\imath\|^{-1}_\mathsf{w}
\]
for any  ${\psi\in\varGamma_\mathsf{w}}$.
In particular, the equality $\varPhi x=\sum x_k\phi_k$  is valid for all  ${x\in E}$. This gives the equalities
\[\|\varPhi x\|_\chi^2=\sum|x_k|^2=\|x\|^2,\quad{x\in E}.\]
Using this, we will examine the composition of  $\varPhi$ with the $\varGamma_\mathsf{w}$-valued function
${\varepsilon\colon E\ni x\mapsto\varepsilon(x)}$.  Its correctness  justifies the following assertion
that substantially uses  the $L^\infty_\chi$-valued function
\[
\phi_x\colon\mathfrak{U}_\pi\ni u\mapsto(\varPhi x)(u)=\sum x_k\phi_k(u)
\]
 which is linear in the variable $x\in E$.

 Similarly to the known case of Wiener spaces,
 the  function $\phi_x$ can be seen as a group analog of the Paley-Wiener map
 (see e.g. \cite[n.4.4]{Sengupta2014} or \cite{Stroock2008}).

\begin{lemma}\label{infty}
The composition $\varPhi\varepsilon(x)$, which is understood as the function
\[
[\varPhi\varepsilon(x)](u)\colon\mathfrak{U}_\pi\ni u\mapsto\exp\left(\phi_x(u)\right),
\]
takes values in  $L_\chi^\infty$ for all ${x\in E}$.
\end{lemma}

\begin{proof}
Applying $\varPhi$ to the Fourier decomposition \eqref{Tayl}, we obtain
\[
\varPhi\varepsilon(x)=\sum_{n\in\mathbb{Z}_+}\frac{1}{n!}
\sum_{(\lambda,\imath)\in\mathbb{Y}_n\times\mathbb{N}^{l(\lambda)}_\ast}
\frac{n!}{\lambda !}{x^\lambda_\imath \phi^\lambda_\imath}
=\sum_{n\in\mathbb{Z}_+}\frac{1}{n!}\Big(\sum_{k\in\mathbb{N}}x_k\phi_k\Big)^n=
\exp\left(\phi_x\right).
\]
It directly follows that $\|\varPhi \varepsilon(x)\|_\infty\le\exp\|\phi_x\|_\infty$.
\end{proof}

\begin{theorem}\label{laplace1}
For every $f={\sum f_n\in{H}^2_\chi}$, $(f_n\in H_\chi^{2,n})$ the entire analytic function
$\widehat{f}(x):={\left\langle\varepsilon(x)\mid\varPhi^* f\right\rangle_\mathsf{w}}$ in the variable ${x\in E}$
and its Taylor coefficients at origin have the integral representations
\begin{equation}\label{laplaceA}
\widehat{f}(x)=\int \exp(\bar\phi_x)f\,d\chi\quad\text{and}\quad d^n_0\widehat{f}(x)=\int\bar\phi_x^nf_n\,d\chi,
\end{equation}
respectively.
The  mapping ${\mathcal{F}\colon{H}^2_\chi\ni{f}\mapsto\widehat{f}\in{H}^2_\mathsf{w}}$
(understanding as a Fock-symmetric $\mathcal{F}$-transform)
 provides the isometries
 \[
 {H^2_\chi\simeq H^2_\mathsf{w}}\quad\text{and}\quad  H_\chi^{2,n}\simeq H^{2,n}_\mathsf{w}.
 \]
\end{theorem}

\begin{proof}
First recall that the $\varGamma_\mathsf{w}$-valued function $\varepsilon(\cdot)$ is entire analytic  on ${E}$,
therefore  $\widehat{f}$ is the same, as the composition of $\varepsilon(\cdot)$ with  ${\langle\cdot\mid{\varPhi^*f}\rangle_\mathsf{w}}$.
Farther on, consider the Fourier decomposition with respect to the basis $\phi^{\mathbb{Y}}$,
\[
f=\sum_{n\in\mathbb{Z}_+}f_n=\sum_{{n\in\mathbb{Z}_+}\atop
{(\lambda,\imath)\in\mathbb{Y}_n\times\mathbb{N}^{l(\lambda)}_\ast}}
\frac{\hat{f}_{\lambda,\imath,n}\bar\phi^\lambda_\imath}{\|\phi^\lambda_\imath\|_\chi},\qquad
\hat{f}_{\lambda,\imath,n}=\frac{1}{\left\|\phi^\lambda_\imath\right\|_\chi}
\int {f\,\bar{\phi}^\lambda_\imath\,d\chi}.
\]
Applying $\varPhi^*$ to $f$ in this decomposition and substituting $\hat{f}_{\lambda,\imath,n}$ into $\widehat{f}$, we obtain
\begin{align*}
\widehat{f}(x)&
=\sum_{n\in\mathbb{Z}_+}\frac{1}{n!}\!\sum_{(\lambda,\imath)\in\mathbb{Y}_n\times\mathbb{N}^{l(\lambda)}_\ast}
\frac{n!}{\lambda !}\frac{\hat{f}_{\lambda,\imath,n}
\langle {e}^{\odot\lambda}_\imath\mid {e}^{\odot\lambda}_\imath\rangle_\mathsf{w}
x^\lambda_\imath}{\|{e}^{\odot\lambda}_\imath\|_\mathsf{w}}\\
&=\int \sum_{n\in\mathbb{Z}_+}\frac{1}{n!}\Big(
\sum_{(\lambda,\imath)\in\mathbb{Y}_n\times\mathbb{N}^{l(\lambda)}_\ast}
\frac{n!}{\lambda !}{x^\lambda_\imath\bar \phi^\lambda_\imath}\Big)f\,d\chi
=\int \exp\left(\bar\phi_x\right)f\,d\chi
\end{align*}
where the last equality is valid by Lemma~\ref{infty}.
It particularly follows that for $y=\alpha{x}$,
\[
\widehat{f}\left(y\right)=
\int \exp\left(\bar\phi_{\alpha x}\right)f\,d\chi=
\sum\alpha^n\!\int \frac{\bar\phi_x^n}{n!}f_nd\chi, \quad{\alpha\in\mathbb{C}}.
\]
Differentiating $\widehat{f}$ at $y=0$ and using the $n$-homogeneity of derivatives, we obtain
\[
d^n_0\widehat{f}(x)=\frac{d^n}{d\alpha^n}\!
\sum\alpha^n\int \frac{\bar\phi_x^n}{n!}f_n\,d\chi
\!\mathrel{\Big|}_{\alpha=0}\!=\int\bar\phi_x^nf_n\,d\chi.
\]
Finally, we notice that the isometry ${H^2_\chi\simeq H^2_\mathsf{w}}$ holds, since
the isometry $\varPhi^*$ is surjective.
In the case of polynomials we similarly get $H_\chi^{2,n}\simeq H^{2,n}_\mathsf{w}$.
\end{proof}

Note that a different integral formula for analytic functions employing Wiener measures on infinite-dimensional Banach
spaces was presented in \cite{PinascoZalduendo05}.

\section{Exponential creation and annihilation groups}\label{5}

Let us define the linear mapping $\mathfrak{j}_n\colon{E}^{\odot n}_\mathsf{w}\rightarrow{E}^{\odot n}_\mathsf{h}$
to be the continuous extension of identity mapping acting on
the dense subspace ${E}^{\odot n}_{\sf alg}\subset{E}^{\odot n}_\mathsf{w}\cap{E}^{\odot n}_\mathsf{h}$.
Such continuous extension $\mathfrak{j}_n$ is a contractive injection with dense range. In fact,
enough to expand elements from ${E}^{\odot n}_\mathsf{w}$ and  ${E}^{\odot n}_\mathsf{h}$
into the Fourier series with respect to orthogonal basis $e^{\odot\mathbb{Y}_n}$ and apply the inequality
\begin{equation}\label{contract}
\|{e}^{\odot\lambda}_\imath\|^2_\mathsf{w}
=C_{|\lambda|,l(\lambda)}\|{e}^{\odot\lambda}_\imath\|_\mathsf{h}^2
\le\|{e}^{\odot\lambda}_\imath\|_\mathsf{h}^2,\quad \lambda\in\mathbb{Y}_n
\end{equation}
which follows from Theorem~\ref{irrep1}, taking into account the inequality
\eqref{constant}. Using subsequently
that ${E}^{\odot n}_\mathsf{h}$ is reflexive, we obtain that its adjoint operator
$\mathfrak{j}_n^*\colon{{E}^{\odot n}_\mathsf{h}\rightarrow{E}^{\odot n}_\mathsf{w}}$ is a contractive
injection with dense range. Thus, the mapping $\mathfrak{j}_n$ is also injective. Moreover,
${E}^{\odot n}_\mathsf{h}\stackrel{\mathfrak{j}_n^*}\rightarrow{E}^{\odot n}_\mathsf{w}
\stackrel{\mathfrak{j}_n}\looparrowright{E}^{\odot n}_\mathsf{h}$ forms a Gelfand triple.
Particularly, the operator $\mathfrak{s}_n$ possesses continuous extension on  ${E}^{\odot n}_\mathsf{w}$.

Using this, we consider the linear operator
\[
\mathfrak{s}_{n/m}:=\mathfrak{s}_n\circ(\mathfrak{j}_m\otimes\mathfrak{j}_{n-m})\quad\text{with}\quad m\le n
\]
 defined to be
 ${\phi_m\odot\psi_{n-m}}={\mathfrak{s}_{n/m}(\phi_m\otimes\psi_{n-m})\in{E}^{\odot n}_\mathsf{w}}$
for all $\phi_m\in{E}^{\odot m}_\mathsf{w}$, $\psi_{n-m}\in{E}^{\odot(n- m)}_\mathsf{w}$.

\begin{lemma}\label{hbar0}
The mapping ${\mathfrak s}_{n/m}$ from
${E}^{\odot m}_\mathsf{w}\otimes_\mathsf{h}{E}^{\odot(n- m)}_\mathsf{w}$ to
${E}^{\odot n}_\mathsf{w}$ is a contractive injection with dense range.
\end{lemma}
\begin{proof}
Expand elements of
${{E}^{\odot m}_\mathsf{w}\otimes_\mathsf{h}{E}^{\odot(n- m)}_\mathsf{w}}$
with respect to  ${e}^{\odot\lambda}_\imath\otimes {e}^{\odot\mu}_\jmath$
for all $\lambda,\mu\in\mathbb{Y}$, ${\imath\in\mathbb{N}^{l(\lambda)}}\!,$ ${\jmath\in\mathbb{N}^{l(\mu)}}$
 such that $|\lambda|=m$, $|\mu|=n-m$. Using \eqref{contract}, we have
\begin{align*}
\|{{e}^{\odot\lambda}_\imath\otimes{e}^{\odot\mu}_\jmath}\|
_{{E}^{\odot m}_\mathsf{w}\otimes_\mathsf{h}{E}^{\odot(n- m)}_\mathsf{w}}&
=\|{e}^{\odot\lambda}_\imath\|_\mathsf{w}\|{e}^{\odot\mu}_\jmath\|_\mathsf{w}\\
&\le\|{e}^{\odot\lambda}_\imath\|_\mathsf{h}\|{e}^{\odot\mu}_\jmath\|_\mathsf{h}
=\|{{e}^{\odot\lambda}_\imath\otimes{e}^{\odot\mu}_\jmath}\|_\mathsf{h}.
\end{align*}
As above,  it implies that the mapping
${\mathfrak{j}_m\otimes\mathfrak{j}_{n-m}}\colon
{{E}^{\odot m}_\mathsf{w}\otimes_\mathsf{h}{E}^{\odot(n- m)}_\mathsf{w}}
\to{E}^{\otimes n}_\mathsf{h}$, defined to be  the continuous extension
of identity mapping on  ${E}^{\odot m}_{\sf alg}\otimes{E}^{\odot(n-m)}_{\sf alg}$,
is a contractive injection. Using subsequently that
${{E}^{\odot m}_\mathsf{h}\otimes_\mathsf{h}{E}^{\odot(n- m)}_\mathsf{h}}$
is reflexive,  we get the Gelfand triple
\[
{{E}^{\odot m}_\mathsf{w}\otimes_\mathsf{h}{E}^{\odot(n- m)}_\mathsf{w}}
\stackrel{\mathfrak{s}_{n/m}}\to{E}^{\odot n}_\mathsf{h}
\stackrel{\mathfrak{j}_n^*}\to{E}^{\odot n}_\mathsf{w}
\]
where injections are contractive and have dense ranges.
\end{proof}

\begin{lemma}\label{expbound0}
The exponential creation group, defined  on  $\left\{\varepsilon(x)\colon x\in{E}\right\}$
by \[\mathcal{T}_\mathsf{a}\varepsilon(x)=\varepsilon(x+\mathsf{a}),\]
has a unique linear extension ${\mathcal{T}_\mathsf{a}\colon
\varGamma_\mathsf{w}\ni\psi\mapsto\mathcal{T}_\mathsf{a}\psi\in\varGamma_\mathsf{w}}$
such that
\[\|\mathcal{T}_\mathsf{a}\psi\|^2_\mathsf{w}\le\exp(\|\mathsf{a}\|^2)\|\psi\|^2_\mathsf{w}\quad
\text{and}\quad\mathcal{T}_\mathsf{a+b}=\mathcal{T}_\mathsf{a}\mathcal{T}_\mathsf{b}=
\mathcal{T}_\mathsf{b}\mathcal{T}_\mathsf{a}\quad\text{for all}\quad{\mathsf{a},\mathsf{b}\in E}.\]
\end{lemma}

\begin{proof}
Let us define the creation operators
$\delta_{\mathsf{a},n}^m\colon{E}^{\odot(n- m)}_\mathsf{w}\rightarrow{E}^{\odot n}_\mathsf{w}$ ${(m\le n)}$  as
\begin{equation}\label{lambda}
\delta_{\mathsf{a},n}^mx^{\otimes(n-m)}:=
\mathfrak{s}_{n/m}\left[\mathsf{a}^{\otimes m}\otimes x^{\otimes (n-m)}\right]
=\frac{(n-m)!}{n!}\frac{d^m(x+t\mathsf{a})^{\otimes n}}{dt^m}{\mathrel{\Big|}}_{t=0}
\end{equation}
for all $\mathsf{a},x\in{E}$.  Note that the second equality in \eqref{lambda}
follows from  the binomial formula for symmetric tensor elements
${(x+t\mathsf{a})^{\otimes n}}=\sum_{m=0}^{n}\binom{n}{m}(t\mathsf{a})^{\otimes m}\odot x^{\otimes (n-m)}$.
Put $\delta_{\mathsf{a},n}^0=1$.  If $\mathsf{a}=0$ then $\delta_{0,n}^m=0$.
Summing over $n\ge m$ with coefficients $1/(n-m)!$, we get
\begin{equation}\label{taylor0}
\delta_\mathsf{a}^m\varepsilon(x)=
\frac{d^m\varepsilon(x+t\mathsf{a})}{dt^m}{\mathrel{\Big|}}_{t=0}
=\bigoplus_{n\ge m}\frac{\mathfrak{s}_{n/m}[\mathsf{a}^{\otimes m}\otimes x^{\otimes (n-m)}]}{(n-m)!},
\quad t\in\mathbb{C}.
\end{equation}
This series is convergent, since by  Lemma~\ref{hbar0} and \eqref{exp} the inequality
\begin{equation}\label{nabla}
\left\|\delta_\mathsf{a}^m\varepsilon(x)\right\|_\mathsf{w}\le\|\mathsf{a}\|^{m}
\Big\|\bigoplus_{n\ge m}\frac{x^{\otimes(n-m)}}{(n-m)!}\Big\|_\mathsf{w}=
\|\mathsf{a}\|^{m}\left\|\varepsilon(x)\right\|_\mathsf{w}
\end{equation}
holds.  From \eqref{taylor0} and the tensor binomial formula  mentioned above it follows that
\begin{align*}
\bigoplus_{m=0}^n\frac{1}{m!}\delta_{\mathsf{a},n}^m\frac{x^{\otimes(n-m)}}{(n-m)!}
=\bigoplus_{m=0}^n\frac{\mathsf{a}^{\otimes m}\odot x^{\otimes (n-m)}}{m!(n-m)!}=\frac{(x+\mathsf{a})^{\otimes n}}{n!}.
\end{align*}
Summing over $n\in\mathbb{Z}_+$ with coefficients $1/n!$ and using \eqref{taylor0}, we obtain
\begin{align*}
\mathcal{T}_\mathsf{a}\varepsilon(x)&=\bigoplus_{n\in\mathbb Z_+}\sum_{m=0}^n\frac{1}{m!}\delta_{\mathsf{a},n}^m\frac{x^{\otimes(n-m)}}{(n-m)!}\\
&=\sum_{m\in\mathbb Z_+}\frac{1}{m!}\bigoplus_{n\ge m}\delta_{\mathsf{a},n}^m\frac{x^{\otimes(n-m)}}{(n-m)!}
=\exp(\delta_\mathsf{a})\varepsilon(x).
\end{align*}
The inequalities \eqref{exp} and  \eqref{nabla} yield
$\left\|\mathcal{T}_\mathsf{a}\varepsilon(x)\right\|^2_\mathsf{w}\le\exp\big(\|\mathsf{a}\|^2\big)
\left\|\varepsilon(x)\right\|^2_\mathsf{w}$.  Taking into account the totality
of $\left\{\varepsilon(x)\colon x\in{E}\right\}$, this inequality
implies the required inequality on $\varGamma_\mathsf{w}$.
It also follows that $\mathcal{T}_\mathsf{a+b}=\mathcal{T}_\mathsf{a}\mathcal{T}_\mathsf{b}=
\mathcal{T}_\mathsf{b}\mathcal{T}_\mathsf{a}$, since
$\delta_{\mathsf{a}+\mathsf{b}}=\delta_\mathsf{a}+\delta_\mathsf{b}$
for all $\mathsf{a},\mathsf{b}\in E$ by linearity of creation operators. This ends the proof.
\end{proof}

We define the adjoint operators
$\delta_{\mathsf{a},n}^{*m}\colon{E}^{\odot n}_\mathsf{w}\ni\psi_n
\rightarrow\delta_{\mathsf{a},n}^{*m}\psi_n\in{E}^{\odot(n- m)}_\mathsf{w}$ as
\begin{equation*}\label{adj0}
\big\langle\delta_{\mathsf{a},n}^mx^{\otimes(n-m)}\mid\psi_n\big\rangle_\mathsf{w}=
\big\langle x^{\otimes (n-m)}\mid\delta_{\mathsf{a},n}^{*m}\psi_n\big\rangle_\mathsf{w},
\quad  \mathsf{a},x\in{E}
\end{equation*}
for ${n\ge m}$.
It immediately  follows that  for every ${\psi_{n-m}\in{E}^{\odot(n- m)}_\mathsf{w}}$ and $x\in{E}$,
\begin{equation}\label{annihilat}\begin{split}
\big\langle\delta_{\mathsf{a},n}^{*m}x^{\otimes n}\mid\psi_{n-m}\big\rangle_\mathsf{w}&=
\big\langle x^{\otimes n}\mid\delta_{\mathsf{a},n}^m\psi_{n-m}\big\rangle_\mathsf{w}
=\big\langle x^{\otimes n}\mid\mathsf{a}^{\otimes m}\odot\psi_{n-m}\big\rangle_\mathsf{w}\\
&=\langle x\mid\mathsf{a}\rangle^m\big\langle x^{\otimes (n-m)}\mid\psi_{n-m} \big\rangle_\mathsf{w}.
\end{split}\end{equation}
Using $\delta_{\mathsf{a},n}^{*m}$, we can define the exponential annihilation group
by the equalities
\begin{equation}\label{star1}
\mathcal{T}_\mathsf{a}^*\varepsilon(x)=\exp(\delta^*_\mathsf{a})\varepsilon(x)=
\sum_{m\in\mathbb Z_+}\frac{\delta_\mathsf{a}^{*m}\varepsilon(x)}{m!},\quad
\delta_\mathsf{a}^{*m}\varepsilon(x):=\bigoplus_{n\ge m}\frac{\delta_{\mathsf{a},n}^{*m}x^{\otimes n}}{n!}
 \end{equation}
 for all $\mathsf{a},x\in{E}$.
Taking into account  Lemma~\ref{expbound0},  we obtain the following claim.

\begin{lemma}\label{lambdadual}
The exponential annihilation group $\mathcal{T}_\mathsf{a}^*$ defined by \eqref{star1} possesses a unique linear extension ${\mathcal{T}^*_\mathsf{a}\colon\varGamma_\mathsf{w}\ni\psi\mapsto\mathcal{T}^*_\mathsf{a}\psi
\in\varGamma_\mathsf{w}}$ such that
 \[
 \left\|\mathcal{T}_\mathsf{a}^*\psi\right\|^2_\mathsf{w}\le
{\exp(\|\mathsf{a}\|^2)}\left\|\psi\right\|^2_\mathsf{w}
\quad\text{and}\quad
\mathcal{T}^*_\mathsf{a+b}=\mathcal{T}^*_\mathsf{a}\mathcal{T}^*_\mathsf{b}=
\mathcal{T}^*_\mathsf{b}\mathcal{T}^*_\mathsf{a}
\quad\text{for all}\quad{\mathsf{a},\mathsf{b}\in E}.\]

\end{lemma}

\section{Intertwining properties of  $\mathcal{F}$-transform}\label{sec:6}

Let us define on  the space $H^2_\chi$  the multiplicative group $M^\dagger_\mathsf{a}\colon{E}\ni\mathsf{a}\mapsto M_\mathsf{a}^\dagger$ to be
\[
M^\dagger_\mathsf{a}f(u)=\exp[\bar\phi_\mathsf{a}(u)] f(u),
\quad  f\in H^2_\chi,\quad u\in\mathfrak U_\pi.
\]
It can be  considered as  a linear representation of the additive group from ${E}$.
By Lemma~\ref{infty} the function $u\mapsto\exp[\bar\phi_\mathsf{a}(u)]$ with a fixed $\mathsf{a}$
belongs to $L^\infty_\chi$. Hence, $M_\mathsf{a}^\dagger$ is continuous on $H^2_\chi$.
The generator of the $1$-parameter group ${\mathbb{C}\ni t\mapsto M^\dagger_{t\mathsf{a}}}$
coincides with the operator of multiplication by the $L_\chi^\infty$-valued function
\[\bar\phi_\mathsf{a}\colon\mathfrak{U}_\pi\ni u\mapsto\bar\phi_\mathsf{a}(u)\quad\text{where}\quad
dM^\dagger_{t\mathsf{a}}/dt|_{t=0}=\bar\phi_\mathsf{a}.\]
The continuity of  ${E}\ni\mathsf{a}\mapsto\exp(\bar\phi_\mathsf{a})$
implies that this $1$-parameter group $M^\dagger_{t\mathsf{a}}$  is strongly continuous on $H^2_\chi$.
As a consequnce, its  generator $(\bar\phi_\mathsf{a} f)(u)=\bar\phi_\mathsf{a}(u) f(u)$ with  domain
$\mathfrak{D}(\bar\phi_\mathsf{a})={\big\{f\in H^2_\chi\colon
\bar\phi_\mathsf{a}f\in  H^2_\chi\big\}}$ is closed and densely-defined.
As well, its  power  $\bar\phi_\mathsf{a}^m$ defined
on $\mathfrak{D}(\bar\phi_\mathsf{a}^m)={\big\{f\in H^2_\chi\colon
\bar\phi_\mathsf{a}^mf\in H^2_\chi\big\}}$ for any ${m\in\mathbb{N}}$
is the same.

The additive group contained in ${E}$ may be also linearly represented on  ${H}^2_\mathsf{w}$ as the shift group
\[
T_\mathsf{a}\widehat{f}(x)=\widehat{f}(x+\mathsf{a}),\quad  f\in H^2_\chi,
\quad x,\mathsf{a}\in{E}.
\]
The directional derivative on the space ${H}^2_\mathsf{w}$
along a nonzero  $\mathsf{a}\in{E}$   coincides with the generator of the $1$-parameter
shift subgroup $\mathbb{C}\ni t\mapsto{T}_{t\mathsf{a}}$, that is,
\[
\mathfrak{d}_\mathsf{a}\widehat{f}=\lim_{t\to 0}t^{-1}(T_{t\mathsf{a}}\widehat{f}-\widehat{f})
\quad\text{with domain}\quad
\mathfrak{D}(\mathfrak{d}_\mathsf{a}):=
\big\{\widehat{f}\in{H}^2_\mathsf{w}\colon\mathfrak{d}_\mathsf{a}\widehat{f}\in{H}^2_\mathsf{w}\big\}.
\]
Note that the $1$-parameter  shift group $T_{t\mathsf{a}}$,
which is intertwined with  $M_{t\mathsf{a}}^\dagger$ by the $\mathcal{F}$-transform
\begin{equation}\label{conec}
T_{t\mathsf{a}}\widehat{f}(x)=\int \exp\left[\bar\phi_{x+t\mathsf{a}}\right]f\,d\chi
=\int \exp(\bar\phi_x)M^\dagger_{t\mathsf{a}}f\,d\chi,
\end{equation}
is strongly continuous on  ${H}^2_\mathsf{w}$.
Since $\mathfrak{D}(\mathfrak{d}_\mathsf{a}^m)$ contains  all polynomials from
${H}^2_\mathsf{w}$, each operator   $\mathfrak{d}_\mathsf{a}^m$ with domain $\mathfrak{D}(\mathfrak{d}_\mathsf{a}^m)={\big\{\widehat{f}\in{H}^2_\mathsf{w}\colon
\mathfrak{d}_\mathsf{a}^m\widehat{f}\in{H}^2_\mathsf{w}\big\}}$
is closed and densely-defined.  From \eqref{conec} it directly follows
\begin{equation}
\begin{split}\label{conec1}
\mathfrak{d}_\mathsf{a}^m\widehat{f}(x)&=
\int \exp(\bar\phi_x)
\frac{d^mM_{t\mathsf{a}}^\dagger}{dt^m}\Big|_{t=0}f\,d\chi
=\int \exp(\bar\phi_x)\bar\phi_\mathsf{a}^mf\,d\chi
\end{split}
\end{equation}
for all  $f\in\mathfrak{D}(\bar\phi_\mathsf{a}^m)$ and $x\in E$.
On the other hand, by Theorem~\ref{laplace1} we have
\begin{equation}\label{T0}
T_\mathsf{a}\widehat{f}(x)=
\left\langle\mathcal{T}_\mathsf{a}\varepsilon(x)\mid\varPhi^*f\right\rangle=
\left\langle\varepsilon(x)\mid\mathcal{T}^*_\mathsf{a}\varPhi^*f\right\rangle=
\int \exp(\bar\phi_x)\varPhi\mathcal{T}^*_\mathsf{a}\varPhi^*f\,d\chi.
\end{equation}
Theorem~\ref{laplace1} together with \eqref{conec} and  \eqref{T0} imply that $M^\dagger_\mathsf{a}$ is connected with the exponential annihilation group $\mathcal{T}^*_\mathsf{a}$
 by the intertwining operator $\varPhi$. This can be written as
 $M^\dagger_\mathsf{a}=\varPhi\mathcal{T}^*_\mathsf{a}\varPhi^*$. Thus,
  the $\mathcal{F}$-transform serves as an intertwining operator for the  groups $M^\dagger_\mathsf{a}$ on  ${H}^2_\chi$. Moreover, using \eqref{conec}, \eqref{conec1} and \eqref{T0}, we obtain
\[
{d^mT_{t\mathsf{a}}\widehat{f}(x)}/{dt^m}|_{t=0}=
\big\langle\varepsilon(x)\mid\delta_\mathsf{a}^{*m}\varPhi^*f\big\rangle_\mathsf{w}
=\mathfrak{d}_\mathsf{a}^m\widehat{f}(x).
\]
As a result, we have proved the following statement.

\begin{theorem}\label{con}
For every $f\in H^2_\chi$ the following equalities hold,
\[T_\mathsf{a}\mathcal{F}(f)=\mathcal{F}(M^\dagger_\mathsf{a}f),\quad
M^\dagger_\mathsf{a}f=\varPhi\mathcal{T}^*_\mathsf{a}\varPhi^*f,
\quad\mathsf{a}\in{E},\]
Moreover,  for every $f\in\mathfrak{D}(\bar\phi_\mathsf{a}^m)$ ${(m\in\mathbb{N})}$
and a nonzero $\mathsf{a}\in{E}$,
\begin{equation*}\label{T1}
\mathfrak{d}_\mathsf{a}^m\widehat{f}(x)=
\left\langle\varepsilon(x)\mid\delta_\mathsf{a}^{*m}\varPhi^*f\right\rangle_\mathsf{w}
=\int \exp(\bar\phi_x)\bar\phi_\mathsf{a}^mf\,d\chi,
\quad x\in{E}.
\end{equation*}
\end{theorem}

Let us consider on ${H}^2_\mathsf{w}$ the multiplicative group with a nonzero  $\mathsf{a}\in{E}$,
\[
M_{\mathsf{a}^*}\widehat{f}(x)=\widehat{f}(x)\exp\langle{x}\mid\mathsf{a}\rangle,
\quad\widehat{f}\in{H}^2_\mathsf{w}.
\]
The generator on ${H}^2_\mathsf{w}$ of the appropriate
$1$-parameter subgroup $\mathbb{C}\ni t\mapsto M_{t\mathsf{a}^*}$  is
\[
{dM_{t\mathsf{a}^*}}/{dt}|_{t=0}=\left\langle\cdot\mid\mathsf{a}\right\rangle
:={\mathsf{a}}^*,\quad \mathsf{a}\in{E}.
\]
Hence,  it coincides with the following linear operator of multiplication
\[({\mathsf{a}}^*\widehat{f})(x)=\left\langle{x}\mid\mathsf{a}\right\rangle\widehat{f}(x)\quad
\text{with domain}\quad\mathfrak{D}({\mathsf{a}}^*)={\big\{\widehat{f}\in{H}^2_\mathsf{w}\colon
{\mathsf{a}}^*\widehat{f}\in{H}^2_\mathsf{w}\big\}}.\]
Its power ${\mathsf{a}}^{*m}$ is densely-defined on
$\mathfrak{D}({\mathsf{a}}^{*m})={\big\{\widehat{f}\in{H}^2_\mathsf{w}\colon
{\mathsf{a}}^{*m}\widehat{f}\in{H}^2_\mathsf{w}\big\}}$
which contains  all polynomials from  ${H}^2_\mathsf{w}$.

Using Lemma~\ref{expbound0} we can represent the additive group from  $E$
over  the space $H^2_\chi$  by the shift group
\[
T_\mathsf{a}^\dagger=\varPhi\mathcal{T}_\mathsf{a}\varPhi^*\quad
\text{with the generator}\quad \delta_\mathsf{a}^\dagger=\varPhi\delta_\mathsf{a}\varPhi^*
\]
defined on $\mathfrak{D}(\delta_\mathsf{a}^\dagger)
=\big\{f\in{H}^2_\chi\colon\delta_\mathsf{a}^{\dagger}{f}\in{H}^2_\chi\big\}$
This means that  $T_\mathsf{a}^\dagger$  is connected via the intertwining operator $\varPhi$
with the exponential creation group $\mathcal{T}_\mathsf{a}$.

\begin{theorem}\label{conn}
For every $f\in H^2_\chi$ the following equality  holds,
\[M_{\mathsf{a}^*}\mathcal{F}(f)=
\mathcal{F}(T_\mathsf{a}^\dagger f),\quad\mathsf{a}\in{E},\] that is,
the $\mathcal{F}$-transform is
 an intertwining operator for  the groups $M_{\mathsf{a}^*}$ on  ${H}^2_\mathsf{w}$ and
 $T_\mathsf{a}^\dagger$ on  ${H}^2_\chi$. Moreover,
for every $f\in\mathfrak{D}(\delta_\mathsf{a}^{\dagger m})={\big\{f\in{H}^2_\chi\colon
\delta_\mathsf{a}^{\dagger m}{f}\in{H}^2_\chi\big\}}$  ${(m\in\mathbb{N})}$
and a nonzero $\mathsf{a}\in{E}$,
\begin{equation}\label{ai}
({\mathsf{a}}^{*m}\widehat{f})(x)
=\left\langle\varepsilon(x)\mid\delta_\mathsf{a}^m\varPhi^*f\right\rangle_\mathsf{w}
=\int \exp(\bar\phi_x)\delta_\mathsf{a}^{\dagger m}f\,d\chi,\quad x\in{E}.
\end{equation}
\end{theorem}

\begin{proof}
The equality \eqref{annihilat} yields $\langle x\mid\mathsf{a}\rangle^m\psi_{n-m}^*(x)=
\big\langle\delta_{\mathsf{a},n}^{*m}x^{\otimes n}\mid\psi_{n-m}\big\rangle_\mathsf{w}$ for all ${n\ge  m}$.
By Theorem~\ref{laplace1} for  any  $f=\sum_nf_n\in{H}^2_\chi$ there exists a unique
$\psi=\bigoplus_n\psi_n$ in $\varGamma_\mathsf{w}$ with $\psi_n\in{E}^{\odot n}$
such that $\varPhi^*f=\psi$ and $f_n=\psi_n^*$.
Summing over all ${m\in\mathbb{Z}_+}$ and ${n\ge m}$  and using \eqref{star1}, we obtain that
\[
\begin{split}
M_{\mathsf{a}^*}\widehat{f}(x)
&=\exp\langle{x}\mid\mathsf{a}\rangle\big\langle\varepsilon(x)\mid\varPhi^*f\big\rangle_\mathsf{w}
=\sum_{m\in\mathbb{Z}_+}\frac{\langle{x}\mid\mathsf{a}\rangle^m}{m!}
\sum_{n\ge m}\psi_{n-m}^*(x)\\
&=\big\langle\mathcal{T}^*_\mathsf{a}\varepsilon(x)\mid\varPhi^*f\big\rangle_\mathsf{w}
=\big\langle\varepsilon(x)\mid\mathcal{T}_\mathsf{a}\varPhi^*f\big\rangle_\mathsf{w}.
\end{split}
\]
By Theorem~\ref{laplace1} and Lemma~\ref{expbound0}  it follows that the equalities
\begin{equation}\label{M0}
M_{t\mathsf{a}^*}\widehat{f}(x)=
\left\langle\varepsilon(x)\mid\mathcal{T}_{t\mathsf{a}}\varPhi^*f\right\rangle_\mathsf{w}
=\int \exp(\bar\phi_x)T_{t\mathsf{a}}^\dagger f\,d\chi,\quad t\in\mathbb{C}
\end{equation}
hold for all $\widehat{f}\in{H}^2_\mathsf{w}$. On the other hand, the equalities \eqref{star1} and \eqref{M0} yield
\begin{align*}
\frac{d^mM_{t\mathsf{a}^*}\widehat{f}(x)}{dt^m}\Big|_{t=0}
=\int \exp(\bar\phi_x)\frac{d^mT_{t\mathsf{a}}^\dagger }{dt^m}\Big|_{t=0}f\,d\chi
=\int \exp(\bar\phi_x)\delta_\mathsf{a}^{\dagger m}f\,d\chi
\end{align*}
for all  $f\in\mathfrak{D}(\delta_\mathsf{a}^{\dagger m})$.
This  in turn yields \eqref{ai}.
\end{proof}

\section{Commutation relations}\label{sec:7}

Describe the commutation relations between
$M_\mathsf{a}^\dagger$  and $T_\mathsf{b}^\dagger$ on the Hardy space ${H}^2_\chi$.

\begin{theorem}\label{comm}
For any nonzero  $\mathsf{a},\mathsf{b}\in E$  the commutation relations
\begin{equation*}\label{comm2}
 M_\mathsf{a}^\dagger T_\mathsf{b}^\dagger=
\exp\langle\mathsf{a}\mid\mathsf{b}\rangle T_\mathsf{b}^\dagger M_\mathsf{a}^\dagger,
\qquad(\bar\phi_\mathsf{a}\delta_\mathsf{b}^\dagger-
\delta_\mathsf{b}^\dagger\bar\phi_\mathsf{a})f={ \langle\mathsf{a}\mid\mathsf{b}\rangle}f
\end{equation*}
 hold, wherein  $f$ belongs to the dense subspace
$\mathfrak{D}(\bar\phi_\mathsf{b}^2)\cap\mathfrak{D}(\delta_\mathsf{a}^{\dagger 2})\subset{H}^2_\chi$.
\end{theorem}

\begin{proof}
Let us prove that  the following equalities hold,
\begin{align}\label{comm1}
T_\mathsf{a}M_{\mathsf{b}^*}&=
\exp\langle\mathsf{a}\mid\mathsf{b}\rangle M_{\mathsf{b}^*}T_\mathsf{a},\qquad
(\mathfrak{d}_\mathsf{a}\mathsf{b}^*-
\mathsf{b}^*\mathfrak{d}_\mathsf{a})\widehat{f}={ \langle\mathsf{a}\mid\mathsf{b}\rangle}\widehat{f}
\end{align}
where  $\widehat{f}\in\mathfrak{D}(\mathsf{b}^{*2})\cap\mathfrak{D}(\mathfrak{d}_\mathsf{a}^2)$.
First property  follows from the direct calculations:
\begin{align*}
M_{\mathsf{b}^*}T_\mathsf{a}\widehat{f}(x)&=
\exp\langle{x}\mid\mathsf{b}\rangle\widehat{f}(x+\mathsf{a}),\\
T_\mathsf{a}M_{\mathsf{b}^*}\widehat{f}(x)&=
\widehat{f}(x+\mathsf{a})\exp\langle{x}\mid\mathsf{b}\rangle
\exp\langle\mathsf{a}\mid\mathsf{b}\rangle
=\exp\langle\mathsf{a}\mid\mathsf{b}\rangle
M_{\mathsf{b}^*}T_\mathsf{a}\widehat{f}(x)
\end{align*}
for all $\widehat{f}\in{H}^2_\mathsf{w}$ and ${x}\in{E}$.
For any
$\widehat{f}\in\mathfrak{D}(\mathsf{b}^{*2})\cap\mathfrak{D}(\mathfrak{d}_\mathsf{a}^2)$
and ${t\in\mathbb{C}}$, we have
\begin{align*}
\frac{d^2}{dt^2}T_{t\mathsf{a}}M_{t\mathsf{b}^*}\widehat{f}\big|_{t=0}&=
\big[\mathfrak{d}^2_\mathsf{a}T_{t\mathsf{a}}M_{t\mathsf{b}^*}\widehat{f}+
2\mathfrak{d}_\mathsf{a}T_{t\mathsf{a}}\mathsf{b}^*M_{t\mathsf{b}^*}\widehat{f}+
T_{t\mathsf{a}}\mathsf{b}^{*2}M_{t\mathsf{b}^*}\widehat{f}\big]_{t=0}\\
&=(\mathfrak{d}_\mathsf{a}^2+2\mathfrak{d}_\mathsf{a}\mathsf{b}^*+\mathsf{b}^{*2})\widehat{f}.
\end{align*}
On the other hand, differentiating again, we have
\begin{align*}
\frac{d}{dt}T_{t\mathsf{a}}M_{t\mathsf{b}^*}\widehat{f}\big|_{t=0}
=\Big[\frac{d}{dt}\exp\langle{t\mathsf{a}}\mid {\bar{t}\mathsf{b}}\rangle
M_{t\mathsf{b}^*}T_{t\mathsf{a}}\widehat{f}+
\exp\langle{t\mathsf{a}}\mid {\bar{t}\mathsf{b}}\rangle
\frac{d}{dt}M_{t\mathsf{b}^*}T_{t\mathsf{a}}\widehat{f}\Big]_{t=0},\\
(\mathfrak{d}_\mathsf{a}^2+2\mathfrak{d}_\mathsf{a}\mathsf{b}^*+\mathsf{b}^{*2})\widehat{f}=
\frac{d}{dt}\Big[\frac{d}{dt}T_{t\mathsf{a}}M_{t\mathsf{b}^*}\widehat{f}\Big]_{t=0}
=\Big[\frac{d^2}{dt^2}\exp\langle{t\mathsf{a}}\mid {\bar{t}\mathsf{b}}\rangle
M_{t\mathsf{b}^*}T_{t\mathsf{a}}\widehat{f}\\
+2\frac{d}{dt}\exp\langle{t\mathsf{a}}\mid {\bar{t}\mathsf{b}}\rangle
\frac{d}{dt}M_{t\mathsf{b}^*}T_{t\mathsf{a}}\widehat{f}
+\exp\langle{t\mathsf{a}}\mid {\bar{t}\mathsf{b}}\rangle
\frac{d^2}{dt^2}M_{t\mathsf{b}^*}T_{t\mathsf{a}}\widehat{f}\Big]_{t=0}\\
=2{ \langle\mathsf{a}\mid\mathsf{b}\rangle}\widehat{f}+
(\mathfrak{d}_\mathsf{a}^2+2\mathsf{b}^*\mathfrak{d}_\mathsf{a}+\mathsf{b}^{*2})\widehat{f}.
\end{align*}
This yields \eqref{comm1}
where  $\mathfrak{D}(\mathsf{b}^{*2})\cap\mathfrak{D}(\mathfrak{d}_\mathsf{a}^2)$ contains
the dense subspace in ${H}^2_\mathsf{w}$  of all polynomials $\widehat{f}$ generating by
finite sums $\varPhi^*(f)=\bigoplus_n\psi_n\in\varGamma_\mathsf{w}$.

Using that $T_\mathsf{b}^\dagger=\mathcal{F}^{-1}M_{\mathsf{b}^*}\mathcal{F}$ and
$M_\mathsf{a}^\dagger =\mathcal{F}^{-1}T_\mathsf{a}\mathcal{F}$ with
$\mathcal{F}^{-1}\colon{H}^2_\mathsf{w}\rightarrow H^2_\chi$
and applying \eqref{comm1}, we obtain
\begin{align*}
&M_\mathsf{a}^\dagger T_\mathsf{b}^\dagger=\mathcal{F}^{-1}T_\mathsf{a}M_{\mathsf{b}^*}\mathcal{F}
=\exp\langle{x}\mid\mathsf{b}\rangle\mathcal{F}^{-1}M_{\mathsf{b}^*}T_\mathsf{a}\mathcal{F}
=\exp\langle{x}\mid\mathsf{b}\rangle T_\mathsf{b}^\dagger M_\mathsf{a}^\dagger,\\
&(\bar\phi_\mathsf{a}\delta_\mathsf{b}^\dagger-\delta_\mathsf{b}^\dagger\bar\phi_\mathsf{a})f=
\mathcal{F}^{-1}(\mathfrak{d}_\mathsf{a}\mathsf{b}^*-\mathsf{b}^*\mathfrak{d}_\mathsf{a})\mathcal{F}f
=\langle\mathsf{a}\mid\mathsf{b}\rangle {f}
\end{align*}
 for all $f\in\mathfrak{D}(\bar\phi_\mathsf{b}^2)\cap\mathfrak{D}(\delta_\mathsf{a}^{\dagger 2})$.
For  any  $f=\sum_nf_n\in{H}^2_\chi$ there exists a unique
$\psi=\bigoplus_n\psi_n$ in $\varGamma_\mathsf{w}$ with $\psi_n\in{E}^{\odot n}_\mathsf{w}$
such that the equalities $\varPhi^*f=\psi$ and $f_n=\psi_n^*$ hold. Hence, the following embedding
${\mathfrak{D}(\bar\phi_\mathsf{b}^2)\cap\mathfrak{D}(\delta_\mathsf{a}^{\dagger 2})}\subset{H}^2_\chi$
is dense.
\end{proof}

\section{Gauss-Weierstrass semigroups}\label{sec:8}

Next we show that the $1$-parame\-ter Gauss-Weierstrass semigroups on the Hardy space ${H}^2_\chi$
can be well described by shift and multiplicative groups (a classic case can be found in
\cite[n.4.3.2]{butzer67}).
For this purpose we use the Gaussian kernel
\[\mathfrak{g}_r(\tau)=\frac{1}{\sqrt{4\pi r}}\exp\Big(-\frac{\tau^2}{4r}\Big),
\quad \tau\in\mathbb{R},\quad r>0.
\]
\begin{theorem}\label{GW}
The $1$-parameter Gauss-Weierstrass semigroups
${\big\{W_r^{\delta_\mathsf{a}^\dagger}\colon r>0\big\}}$  and
${\big\{W_r^{\bar\phi_\mathsf{a}}\colon r>0\big\}}$, defined on the Hardy space
${H}^2_\chi$ for any nonzero $\mathsf{a}\in{E}$ as
\begin{equation}\label{M2}
W_r^{\delta_\mathsf{a}^\dagger}f=
\int_\mathbb{R}\mathfrak{g}_r(\tau)T_{\tau\mathsf{a}}^\dagger{f}\,d\tau
\quad\text{and}\quad
W_r^{\bar\phi_\mathsf{a}}f=
\int_\mathbb{R}\mathfrak{g}_r(\tau)M_{\tau\mathsf{a}}^\dagger{f}\,d\tau,\quad {f}\in{H}^2_\chi,
\end{equation}
are generated by $\delta_\mathsf{a}^{\dagger 2}$ and $\bar\phi_\mathsf{a}^2$, respectively.
\end{theorem}

\begin{proof}
First it is sufficient to prove that the axillary $1$-parameter families of linear operators over  ${H}^2_\mathsf{w}$
\begin{equation}\label{M}
G_r^{\mathsf{a}^*}\widehat{f}=
\int_\mathbb{R}\mathfrak{g}_r(\tau)M_{\tau\mathsf{a}^*}\widehat{f}\,d\tau
\quad\text{and}\quad
G_r^{\partial_\mathsf{a}}\widehat{f}=
\int_\mathbb{R}\mathfrak{g}_r(\tau)T_{\tau\mathsf{a}}\widehat{f}\,d\tau,\quad \widehat{f}\in{H}^2_\mathsf{w}
\end{equation}
can be generated by $\mathsf{a}^{*2}$ and $\mathfrak{d}_\mathsf{a}^2$ and
satisfy the semigroup property. Properties of Gaussian kernel  yield
\[\begin{split}
\int_\mathbb{R}\mathfrak{g}_r(\tau)\tau^{2k}\,d\tau&=
\frac{1}{2\sqrt{\pi r}}\int_\mathbb{R}
e^{-\frac{\tau^2}{4r}}\tau^{2k}d\tau\Big|_{\tau=2\sqrt{r}\upsilon}
=\frac{(2\sqrt{r})^{2k}}{\sqrt{\pi}}\int_\mathbb{R}e^{-\upsilon^2}\upsilon^{2k}\,d\upsilon\\
&=\frac{2^{2k}r^k}{\sqrt{\pi}}\Gamma\left(\frac{2k+1}{2}\right)=\frac{2(2k-1)!}{(k-1)!}r^k,\quad
k\in\mathbb N.
\end{split}\]
We can rewrite $G^{\mathsf{a}^*}_r\widehat{f}$
on the dense subspace $\big\{\widehat{f}\in{H}^2_\mathsf{w}\colon
\exp(\tau\mathsf{a}^*)\widehat{f}\in{H}^2_\mathsf{w}\big\}$ as
\begin{align*}\label{a}
G^{\mathsf{a}^*}_r\widehat{f}&
=\int_\mathbb{R}\mathfrak{g}_r(\tau)\exp(\tau\mathsf{a}^*)\widehat{f}\,d\tau
=\sum_{l\in\mathbb Z_+}\frac{\mathsf{a}^{*l}\widehat{f}}{l!}
\int_\mathbb{R}\mathfrak{g}_r(\tau)\tau^l\,d\tau\\\nonumber
&=\sum_{k\in\mathbb{Z}_+}\dfrac{2(2k-1)!}{(k-1)!}\frac{r^k\mathsf{a}^{*2k}\widehat{f}}{(2k)!}
=\sum_{k\in\mathbb Z_+}\frac{r^k\mathsf{a}^{*2k}\widehat{f}}{k!}
=\exp(r\mathsf{a}^{*2})\widehat{f}
\end{align*}
By first equality in \eqref{M} the family $G^{\mathsf{a}^*}_r$ can be extended to the convolution
 \[\mathfrak{g}_r\circledast\widehat{f}:=\int_\mathbb{R}\mathfrak{g}_r(\tau)M_{\tau\mathsf{a}^*}\widehat{f}\,d\tau,
\quad \widehat{f}\in{H}^2_\mathsf{w}\]
(dependent on  $\mathsf{a}$) over the whole space ${H}^2_\mathsf{w}$.
Thus, to show that the semigroup property holds, it suffices to show that
\[
\mathfrak{g}_{r+s}\circledast\widehat{f}=G^{\mathsf{a}^*}_{r+s}\widehat{f}
=(G^{\mathsf{a}^*}_r\circ G^{\mathsf{a}^*}_{s})\widehat{f}
=\mathfrak{g}_r\circledast(\mathfrak{g}_s\circledast\widehat{f})=(\mathfrak{g}_r*\mathfrak{g}_s)\circledast\widehat{f}.
\]
But this straightly follows from the known  convolution equality
$\mathfrak{g}_{r+s}=\mathfrak{g}_r*\mathfrak{g}_s$.

Further, using the equality $T_\mathsf{a}^\dagger=\mathcal{F}^{-1}M_{\mathsf{a}^*}\mathcal{F}$ we obtain that
\begin{align*}\label{a}
W_r^{\delta_\mathsf{a}^\dagger}f&
=\int_\mathbb{R}\mathfrak{g}_r(\tau)\mathcal{F}^{-1}M_{\tau\mathsf{a}^*}\mathcal{F}{f}\,d\tau=
\mathcal{F}^{-1}G_r^{\mathsf{a}^*}\mathcal{F}f
\end{align*}
 for all $f\in H^2_\chi$. By Theorem~\ref{conn} it follows that
\[
\frac{dW_r^{\delta_\mathsf{a}^\dagger}f}{dr}\Big|_{r=0}=
\mathcal{F}^{-1}\frac{G_r^{\mathsf{a}^*}\widehat{f}}{dr}\Big|_{r=0}=
\mathcal{F}^{-1}\mathsf{a}^{*2}\widehat{f}=\delta_\mathsf{a}^{\dagger 2}f
\]
for all $f\in\mathfrak{D}(\delta_\mathsf{a}^{\dagger 2})$, since
$\widehat{f}\in\mathfrak{D}(\mathsf{a}^{*2})$ and
$\delta_\mathsf{a}^{\dagger 2}=\mathcal{F}^{-1}\mathsf{a}^{*2}\mathcal{F}$.
Hence, the case of semigroup $W_r^{\delta_\mathsf{a}^\dagger}$ is proven.

Similar reasonings can be applied  to the semigroup $G_r^{\partial_\mathsf{a}}$.
As a result, we obtain that the equalities
$W_r^{\bar\phi_\mathsf{a}}=\mathcal{F}^{-1}G_r^{\partial_\mathsf{a}}\mathcal{F}$
and $\bar\phi_\mathsf{a}^2=\mathcal{F}^{-1}{\mathfrak{d}_\mathsf{a}^2}\mathcal{F}$ hold.
\end{proof}

\section{Complexified infinite-dimensional Heisenberg group}\label{sec:9}

Let us give yet another application.
Consider an infinite-dimensional analog of the Heisenberg group
over $\mathbb{C}$.  Namely, let us define  the group $\mathfrak{G}$ of upper
triangular matrix-type elements
\[X(\mathsf{a},\mathsf{b},t)=\begin{bmatrix}
                      1 & \mathsf{a}&t \\
                      0 & 1&\mathsf{b} \\
                      0&0&1
                    \end{bmatrix},\quad t\in\mathbb{C},\quad \mathsf{a},\mathsf{b}\in{E}
\]
with unit $X(0,0,0)$ and multiplication
\[
\begin{bmatrix}
                      1 & \mathsf{a}&t \\
                      0 & 1&\mathsf{b} \\
                      0&0&1
                    \end{bmatrix}
  \begin{bmatrix}
                      1 & \mathsf{a}'&t' \\
                      0 & 1&\mathsf{b}' \\
                      0&0&1
                    \end{bmatrix}=
 \begin{bmatrix}
                      1 &\mathsf{a}+\mathsf{a}'&t+t'+\langle\mathsf{a}\mid\mathsf{b}'\rangle\\
                      0 & 1&\mathsf{b}+ \mathsf{b}'\\
                      0&0&1
                    \end{bmatrix}.
\]
Obviously, $X(\mathsf{a},\mathsf{b},t)^{-1}=X(-\mathsf{a},-\mathsf{b},-t+\langle\mathsf{a}\mid\mathsf{b}\rangle)$.

Describe an irreducible linear representation of the group $\mathfrak{G}$.
For this purpose we will use the algebra $\mathbb{H}$ of quaternions
$
\gamma={\alpha_1+\alpha_2\mathbbm{i}+\beta_1\mathbbm{j}+\beta_2\mathbbm{k}}
={{(\alpha_1+\alpha_2\mathbbm{i})}+{(\beta_1+\beta_2\mathbbm{i})}\mathbbm{j}}
={\alpha+\beta\mathbbm{j}}
$
as pairs of complex numbers $(\alpha,\beta)\in\mathbb{C}^2$ with
$\alpha={\alpha_1+\alpha_2\mathbbm{i}},\beta={\beta_1+\beta_2\mathbbm{i}\in\mathbb{C}}$ and
$\alpha_\imath,\beta_\imath\in\mathbb{R}$ $(\imath=1,2)$
where basis elements in $\mathbb{R}^4$ satisfy the  relations
$\mathbbm{i}^2=\mathbbm{j}^2=\mathbbm{k}^2=
\mathbbm{i}\mathbbm{j}\mathbbm{k}=-1$,
$\mathbbm{k}=\mathbbm{i}\mathbbm{j}=-\mathbbm{j}\mathbbm{i}$,
$\mathbbm{k}\mathbbm{i} = −\mathbbm{i}\mathbbm{k} = \mathbbm{j}$.
Thus,  $\mathbb{H}=\mathbb{C}\oplus\mathbb{C}\mathbbm{j}$ is a vector space over
$\mathbb{C}$  \cite{zbMATH01028010}.
Denote $\beta:=\Im{\gamma}$ where $\gamma=\alpha+\beta\mathbbm{j}$.

Let $E_\mathbb{H}=E\oplus E\mathbbm{j}$ be the Hilbert  space
 with   $\mathbb{H}$-valued scalar product
\begin{align*}
&\langle\mathsf{p}\mid\mathsf{p}'\rangle=
\langle\mathsf{a}+\mathsf{b}\mathbbm{j}\mid
\mathsf{a}'+\mathsf{b}'\mathbbm{j}\rangle
=\langle\mathsf{a}\mid\mathsf{a}'\rangle+
\langle\mathsf{b}\mid\mathsf{b}'\rangle+
\left[\langle\mathsf{a}'\mid\mathsf{b}\rangle-
\langle\mathsf{a}\mid\mathsf{b}'\rangle\right]\mathbbm{j}
\end{align*}
where  $\mathsf{p}=\mathsf{a}+\mathsf{b}\mathbbm{j}$ with
$\mathsf{a}$, $\mathsf{b}\in{E}$ (similarly, for $\mathsf{p}'=\mathsf{a}'+\mathsf{b}'\mathbbm{j}$).
Hence,
\[
\Im\langle\mathsf{p}\mid\mathsf{p}'\rangle=
\langle\mathsf{a}'\mid\mathsf{b}\rangle-
\langle\mathsf{a}\mid\mathsf{b}'\rangle,\qquad
\Im\langle\mathsf{p}\mid\mathsf{p}\rangle=0.
\]

\begin{theorem}\label{heis}
The following linear representation of $\mathfrak{G}$ over ${H}^2_\chi$
(which can be seen as an analog of the  Weyl-Schr\"{o}dinger representation),
\[
\mathscr{W}^\dagger\colon\mathfrak{G}\ni X(\mathsf{a},\mathsf{b},t)\longmapsto
\exp\Big[t+\frac{1}{2}\langle\mathsf{a}\mid\mathsf{b}\rangle\Big]
 T_\mathsf{a}^\dagger M_\mathsf{b}^\dagger,
\]
is well defined and irreducible.
\end{theorem}

\begin{proof}
First we prove that the following operator representation
\[
\mathscr{W}\colon\mathfrak{G}\ni X(\mathsf{a},\mathsf{b},t)\longmapsto
\exp\Big[t+\frac{1}{2}\langle\mathsf{a}\mid\mathsf{b}\rangle\Big]
M_{\mathsf{a}^*}T_{\mathsf{b}}
\]
into the operator algebras  over ${H}^2_\mathsf{w}$ is well defined and irreducible.
Consider the auxiliary group $\mathbb{C}\times E_\mathbb{H}$ with the multiplication
\begin{align*}
(t,\mathsf{p})(t',\mathsf{p}')&=
\Big(t+t'-\frac{1}{2}\Im\langle\mathsf{p}\mid\mathsf{p}'\rangle,\,
\mathsf{p}+\mathsf{p}'\Big)
\end{align*}
for all $\mathsf{p}=\mathsf{a}+\mathsf{b}\mathbbm{j}$,
$\mathsf{p}'=\mathsf{a}'+\mathsf{b}'\mathbbm{j}\in E_\mathbb{H}$.
It  is related to $\mathfrak{G}$ via the mapping
\[
\mathscr{G}\colon
X(\mathsf{a},\mathsf{b},t)\longmapsto
\Big(t-\frac{1}{2}\langle\mathsf{a}\mid\mathsf{b}\rangle, \,
\mathsf{a}+\mathsf{b}\mathbbm{j}\Big).
\]
Check that $\mathscr{G}$ is a group isomorphism. In fact,
\begin{align*}
&\mathscr{G}\left(X(\mathsf{a},\mathsf{b},t)X(\mathsf{a}',\mathsf{b}',t')\right)
=\mathscr{G}\left(X(\mathsf{a}+\mathsf{a}',\mathsf{b}+\mathsf{b}',
t+t'+\langle\mathsf{a}\mid\mathsf{b}'\rangle)\right)\\
&=\Big(t+t'+\langle\mathsf{a}\mid\mathsf{b}'\rangle-\frac{1}{2}\left[
\langle\mathsf{a}+\mathsf{a}'\mid\mathsf{b}+\mathsf{b}'\rangle\right],(\mathsf{a}+\mathsf{a}')+
(\mathsf{b}+\mathsf{b}')\mathbbm{j}\Big)\\
&=\Big(t+t'-\frac{1}{2}\big[\langle\mathsf{a}\mid\mathsf{b}\rangle+
\langle\mathsf{a}'\mid\mathsf{b}'\rangle\big]
+\frac{1}{2}\big[\langle\mathsf{a}\mid\mathsf{b}'\rangle-
\langle\mathsf{a}'\mid\mathsf{b}\rangle\big],(\mathsf{a}+\mathsf{a})+
(\mathsf{b}+\mathsf{b}')\mathbbm{j}\Big)\\
&=\Big(t-\frac{1}{2}\langle\mathsf{a}\mid\mathsf{b}\rangle, \,
\mathsf{a}+\mathsf{b}\mathbbm{j}\Big)
\Big(t'-\frac{1}{2}\langle\mathsf{a}'\mid\mathsf{b}'\rangle, \,
\mathsf{a}'+\mathsf{b}'\mathbbm{j}\Big)\\
&=\mathscr{G}\left(X(\mathsf{a},\mathsf{b},t)\right)
\mathscr{G}\left(X(\mathsf{a}',\mathsf{b}',t')\right).
\end{align*}

Now let us check that  the Weyl-like operator
\[
W(\mathsf{p})=
\exp\Big[\frac{1}{2}\langle\mathsf{a}\mid\mathsf{b}\rangle\Big]M_{\mathsf{a}^*}T_{\mathsf{b}},
\quad\mathsf{p}=\mathsf{a}+\mathsf{b}\mathbbm{j}
\]
on the space ${H}^2_\mathsf{w}$ satisfies the commutation relation
\[
W(\mathsf{p}+\mathsf{p}')=
\exp\Big[-\frac{1}{2}\Im\left\langle\mathsf{p}\mid\mathsf{p}'\right\rangle\Big]
W(\mathsf{p})W(\mathsf{p}').
\]
In fact, using \eqref{comm1}, we obtain
\begin{align*}
&\exp\Big[\frac{1}{2}\langle\mathsf{a}\mid\mathsf{b}'\rangle
-\frac{1}{2}\langle\mathsf{a}'\mid\mathsf{b}\rangle\Big]W(\mathsf{p})W(\mathsf{p}')\\
&=\exp\Big[\frac{1}{2}\langle\mathsf{a}\mid\mathsf{b}\rangle
+\frac{1}{2}\langle\mathsf{a}'\mid\mathsf{b}'\rangle\Big]
\exp\Big[\frac{1}{2}\langle\mathsf{a}\mid\mathsf{b}'\rangle
-\frac{1}{2}\langle\mathsf{a}'\mid\mathsf{b}\rangle\Big]
M_{\mathsf{a}^*}T_\mathsf{b}M_{{\mathsf{a}'}^*}T_{\mathsf{b}'}\\
&=\exp\Big[\frac{1}{2}\langle\mathsf{a}+\mathsf{a}'\mid\mathsf{b}+\mathsf{b}'\rangle\Big]
M_{\mathsf{a}^*+{\mathsf{a}'}^*}T_{\mathsf{b}+\mathsf{b}'}
=W(\mathsf{p}+\mathsf{p}').
\end{align*}
As a consequence,  the  mapping
${\mathscr{I}\colon\mathbb{C}\times E_\mathbb{H}\ni(t,\mathsf{p})\longmapsto\exp(t)W(\mathsf{p})}$
is a group isomorphism. So, $\mathscr{W}$ is also a group isomorphism as a composition
of  the group isomorphisms $\mathscr{I}$ and  $\mathscr{G}$.

Let us check irreducibility. If there exists an element $x_0\noteq 0$ in $E$ and an integer ${n>0}$ such that
\[
\exp\Big[t+\frac{1}{2}\langle\mathsf{a}\mid\mathsf{b}\rangle\Big]
e^{\langle\mathsf{a}\mid x\rangle}\left[x_0^*(x+\mathsf{b})\right]^n=0\quad\text{for all}\quad
x,\mathsf{a},\mathsf{b}\in E
\]
then $x_0=0$. This gives a contradiction.
Hence the representation $\mathscr{W}$ is irreducible. Finally, using that
\[
\exp\Big[t+\frac{1}{2}\langle\mathsf{a}\mid\mathsf{b}\rangle\Big]
 T_\mathsf{a}^\dagger M_\mathsf{b}^\dagger=
 \mathcal{F}^{-1}\Big(\exp\Big[t+\frac{1}{2}\langle\mathsf{a}\mid\mathsf{b}\rangle\Big]
M_{\mathsf{a}^*}T_{\mathsf{b}}\Big)\mathcal{F},
\]
we conclude that the group representation $\mathscr{W}^\dagger=\mathcal{F}^{-1}\mathscr{W}\mathcal{F}$ is  irreducible.
\end{proof}


\bibliographystyle{spmpsci}      

\end{document}